\documentclass{amsart}
\usepackage{amsthm,amsfonts,amsmath,amscd,amssymb,latexsym,epsfig}

\newcommand{\cx}{{\mathbb C}}

\newcommand{\ad}{\operatorname{ad}}

\newcommand{\im}{\operatorname{Im}}

\newcommand{\End}{\operatorname{End}}

\newcommand{\Ker}{\operatorname{Ker}}

\newcommand{\wh}{\widehat}
\newcommand{\ol}{\overline}

\numberwithin{equation}{section}

\newtheorem{theorem}{Theorem}[section]

\newtheorem{proposition}[theorem]{Proposition}

\theoremstyle{remark}

\newtheorem{remark}[theorem]{Remark}

\newtheorem{definition}[theorem]{Definition}

\newcommand{\oN}{{\mathbb{N}}}

\newcommand{\oP}{{\mathbb{P}}}

\newcommand{\oR}{{\mathbb{R}}}

\newcommand{\oZ}{{\mathbb{Z}}}

\newcommand{\sC}{{\mathcal{C}}}   
\newcommand{\sD}{{\mathcal{D}}}

\newcommand{\sF}{{\mathcal{F}}}
\newcommand{\sH}{{\mathcal{H}}}
\newcommand{\sI}{{\mathcal{I}}}
\newcommand{\sJ}{{\mathcal{J}}}

\newcommand{\sL}{{\mathcal{L}}}   
\newcommand{\sM}{{\mathcal{M}}}   
\newcommand{\sN}{{\mathcal{N}}}
\newcommand{\sO}{{\mathcal{O}}}

\newcommand{\sS}{{\mathcal{S}}}

\newcommand{\sV}{{\mathcal{V}}}

\newcommand{\sX}{{\mathcal{X}}}

\newcommand{\sZ}{{\mathcal{Z}}}

\begin{document}

\title[Limits of monopoles and of pluricomplex structures]{Hypercomplex limits of pluricomplex structures and the Euclidean limit of hyperbolic monopoles}
\author{Roger Bielawski \& Lorenz Schwachh\"ofer}
\address{School of Mathematics\\
University of Leeds\\Leeds LS2 9JT\\ UK}
\address{Fakult\"at f\"ur Mathematik\\ TU Dortmund\\
	D-44221 Dortmund\\ Germany}



\begin{abstract} We discuss the Euclidean limit of hyperbolic $SU(2)$-monopoles, framed at infinity, from the point of view of pluricomplex geometry.
More generally, we discuss the geometry of hypercomplex manifolds arising as limits of pluricomplex manifolds.
\end{abstract}

\maketitle

\thispagestyle{empty}

\section{Introduction}

BPS monopoles have especially nice moduli spaces on the Euclidean and the hyperbolic $3$-spaces. In particular,  moduli spaces of framed Euclidean monopoles are complete hyperk\"ahler manifolds (at least for the gauge group $SU(2)$, or $SU(N)$ and maximal symmetry breaking at infinity). On the other hand, the geometry of hyperbolic monopoles has been identified in \cite{BS} (building on an earlier work by Nash \cite{Nash}), as  {\em pluricomplex geometry}, i.e. given by a $2$-sphere of complex structures without any anticommuting properties. In the present work, we aim to show how the pluricomplex geometry of hyperbolic monopoles converges to the hyperk\"ahler geometry of Euclidean monopoles, as the curvature of the hyperbolic space tends to zero. An unexpected result of our analysis is that the geometry of Euclidean monopoles is richer than hyperk\"ahler geometry: in addition to a $2$-sphere of complex structures, there is, for each of these complex structures $J$, a preferred endomorphism $A_J$ of the tangent bundle, which anticommutes with $J$ (and usually is not a complex structure). These endomorphisms $A_J$ satisfy integrability conditions, and their existence can be viewed as a purely differential-geometric explanation of why the hypercomplex geometry of Euclidean monopoles is determined by algebraic curves of higher genus.
\par
The existence of such endomorphisms $A_J$ is a more general phenomenon: they arise for any hypercomplex structure which is a limit of pluricomplex structures. We call such structures {\em l-hypercomplex} (``l'' for limit), and discuss them in Section \ref{LHC}. We think that they will repay further investigation, particularly those which we call {\em  strongly integrable}, i.e. those, for which the endomorphisms $A_J$ satisfy integrability conditions (just as for moduli spaces of Euclidean monopoles).
\par
We finish this introduction with some bibliographical comments.
The Euclidean limit of hyperbolic monopoles has been discussed in detail by Jarvis and Norbury \cite{JN}, albeit with a different framing (in the interior  rather than on the boundary). With our choice of framing (which is the one that leads to hyperk\"ahler or pluricomplex geometry), the twistorial description of the Euclidean limit has been given by Atiyah \cite{A}. Finally, the limiting behaviour of spectral curves for charge $2$ monopoles has been described by Norbury and Rom\~ao \cite{NR}.

\section{Convergence of minitwistor spaces}

Hitchin \cite{Hit1} has shown that the space of oriented geodesics in a $3$-manifold of constant curvature has a natural structure of a complex manifold. For $\oR^3$, this complex manifold is $T\oP^1$, while for the hyperbolic space $H^3$ it is $\oP^1\times \ol{\oP^1}-\Delta$, where the two copies of $\oP^1$ correspond to the startpoint and the endpoint of a geodesic on the ideal boundary $S^2\simeq \oP^1$  of $H^3$.
\par
 We find it convenient, as do many other authors, to change the complex structure on the second factor to the opposite one, using the antipodal map. Thus, we identify the space of geodesics in $H^3$ with $\oP^1\times \oP^1-\ol\Delta$, where $\ol\Delta$ is the antidiagonal.
\par
There is a subtlety here: the identification of $\ol{\oP^1}$ with $\oP^1$ implies  a choice of a preferred point in $H^3$. Indeed,  we could use, instead of the antipodal map, an arbitrary antiholomorphic M\"obius transformation  $f:\oP^1\to  \ol{\oP^1}$. The geodesics beginning at $\zeta\in \oP^1$ and ending at $f(\zeta)$ have a unique common point $x_f$. Observe that, under the identification of points in $H^3$ with real curves in $\oP^1\times \oP^1- \ol{\Delta}_f$ ($\ol{\Delta}_f=\{\zeta,f(\zeta);\zeta\in \oP^1\}$), the point $x_f$ always corresponds to the diagonal in $\oP^1\times \oP^1$.
\par
In the case when $f$ is the antipodal map, we denote $x_f$ by $x_0$. We now rescale the hyperbolic space with $x_0$ as preferred point. In other words,
we consider the family $(H^3_\kappa,x_0)$ of pointed hyperbolic spaces, with the scalar curvature of $H^3_\kappa$ equal to $-\kappa$. Its limit, as $\kappa\to 0$, is the Euclidean $3$-space $(\oR^3,0)$. We want to describe the limiting behaviour of the spaces of geodesics as complex manifolds. Observe, that if $\gamma_n$ is a sequence of geodesics in $ H^3_{\kappa_n}$, corresponding to $(\zeta_n,\eta_n)\in \oP^1\times \oP^1-\ol\Delta$, and converging to a Euclidean geodesic $\gamma_\infty$, then $(\zeta_n,\gamma_n)$ convergences to a point $(\zeta_\infty,\zeta_\infty)$ on the diagonal, since the two asymptotic directions of a Euclidean geodesic are opposite to each other. It is not hard to see that the limiting Euclidean geodesic $\gamma_\infty$ is the one we expect: $\gamma_\infty$ corresponds to the point $v_\infty \in T_{\zeta_\infty}\Delta$, where $v_\infty$ is the limit of the direction between $\zeta_n$ and $\eta_n$ (this has been already observed by Norbury and Rom\~ ao \cite[p.310]{NR}). In other words, if we choose coordinates $(\zeta,w)$ on $T\oP^1$ so that $\zeta$ is the affine coordinate on $\oP^1$ and $w\frac{\partial}{\partial \zeta}$ are the corresponding vector fields, then the geodesic $\gamma_\infty$ has coordinates $(\zeta_\infty, w_\infty)$, where $w_\infty=\lim\frac{\eta_n-\zeta_n}{\sqrt{\kappa_n}}$. The reason for the square root is that the rescaling by $t$ changes the scalar curvature by $t^2$.
\par
The convergence of minitwistor spaces has been also described (in slightly different terms) in Atiyah's paper \cite[p.19--21]{A}.

\subsection{A holomorphic description\label{sX}} We now give a description of the above in terms of deformation theory. We shall construct  a complex $3$-manifold $\sX$ fibreing over $\cx$, so that the fibre over $0$ is $T\oP^1$, all other fibres are isomorphic to $ \oP^1\times \oP^1-\ol\Delta$. More precisely, we require $\sX$ and the holomorphic projection 
 $\pi: \sX\to \cx$ to have the following properties:
\begin{itemize}
\item $\pi$ is a holomorphic submersion.
 \item $\pi^{-1}(0)$ is isomorphic to $T\oP^1$.
\item There exists a biholomorphism $\phi:\sX-\pi^{-1}(0)\to \bigl(\oP^1\times \oP^1-\ol\Delta\bigr)\times \cx^\ast$, the restriction $\phi_t$ of which to each fibre $\pi^{-1}(t)$ maps this fibre to $\bigl(\oP^1\times \oP^1-\ol\Delta\bigr)\times \{t\}$.
\item If  $x_n$ is a convergent sequence of points in $\sX$, such that $t_n=\pi(x_n)\to 0$, then $\phi_{t_n}(x_n)$ converges to a point on the diagonal $\Delta\subset \oP^1\times \oP^1$.
\item There is a real structure on $\sX$, compatible with the above.
\end{itemize}

\begin{remark} The last property means that the real structure $\sigma$ on $\sX$ induces the standard real structure on $\pi^{-1}(0)\simeq T\oP^1$, $\pi$ intertwines $\sigma$ and the complex conjugation on $\cx$, and $\phi$ intertwines $\sigma$ and the standard real structure on $\bigl(\oP^1\times \oP^1-\ol\Delta\bigr)\times \cx^\ast$ given by 
\begin{equation}
 \sigma_0(\zeta,\eta,t)=(-1/\bar\eta,-1/\bar\zeta,\bar t).\label{si0}
\end{equation}
\end{remark}

\medskip

We construct $\sX$ by glueing two charts 
$$ U_0=\{(\zeta,w,t)\in \cx^3;\, |\zeta|^2+t w \bar\zeta+1\neq 0\},$$
$$ U_1=\{(\tilde\zeta,\tilde w,\tilde t)\in \cx^3;\, |\tilde\zeta|^2+\tilde t\tilde w \bar{\tilde\zeta}+1\neq 0\},$$
over $\zeta,\tilde\zeta\neq 0$, $\zeta+tw, \tilde\zeta+\tilde t\tilde w\neq 0$ via the transition map:
\begin{equation}
 \tilde\zeta=\frac{1}{\zeta},\quad \tilde w=\frac{-w}{(\zeta+ tw)\zeta},\quad \tilde{t}=t. \label{transition}
\end{equation}

We have a holomorphic projection $\pi:\sX\to \cx$, $\pi(\zeta,w,t)=t$ and we write $\sX_t$ for $\pi^{-1}(t)$. The fibre $\sX_0$ is $T\oP^1$. For $t\neq 0$ we can change 
the coordinates $(\zeta,w)$ and $(\tilde\zeta,\tilde w)$ to  $(\zeta,\eta)$ and $(\tilde\zeta,\tilde \eta)$, where
\begin{equation}
 \eta=\zeta+t w,\enskip \tilde{\eta}=\tilde{\zeta}+t\tilde w. \label{eta}
\end{equation}
It follows then that $\tilde\eta=\eta^{-1}$ and, consequently $\sX_t$ is biholomorphic to  $\oP^1\times \oP^1 -\ol\Delta$ for $t\neq 0$. The map \eqref{eta} describes the biholomorphism $\phi$:
\begin{equation}
 \phi(\zeta,w,t)=(\zeta,\eta,t).\label{phi}
\end{equation}

Finally, the real structure on $\sX$ is defined by
\begin{equation}
 \sigma(\zeta,w,t)=\left(\frac{-1}{\bar\zeta+\bar t\bar{w}}\,,\,\frac{-\bar{w}}{(\bar\zeta+\bar t\bar{w})\bar\zeta}\,,\, \bar t\right).\label{real}
\end{equation}
In coordinates $(\zeta,\eta,t)$, $t\neq 0$, $\sigma$ is given by $\sigma(\zeta,\eta,t)=(-1/\bar\eta,-1/\bar\zeta,\bar t)$, and so it coincides with $\sigma_0$, given by \eqref{si0}.

\begin{remark} We can compactify the $t$-variable by adding $\oP^1\times \oP^1-\ol\Delta$ corresponding to $t=\infty$. We denote the resulting complex manifold by $\bar\sX$. The fibration $\pi$ extend to $\pi:\bar\sX\to \oP^1$, and the real structure $\sigma$ extends to $\bar\sX$.\label{cpct}\end{remark} 

\begin{remark} The manifold $\sX$  can be realised as an open subset of $\oP^3$ with the induced real structure and it occurs in that guise in the paper of Atiyah \cite{A}. Our description, given above, is more suitable to pluricomplex geometry.\end{remark}

\subsection{Convergence}
We wish to make explicit the relation between curves in $\sX$ and points in hyperbolic spaces of varying curvature. We identify $H^3_\kappa$ (the hyperbolic space of scalar curvature $-\kappa$) with the half space $\{(x,y,z)\in \oR^3;\, z>- \sqrt{-1/\kappa}\}$, equipped with the metric
\begin{equation} ds^2=\frac{r^2}{(z+r)^2}(dx^2+dy^2+dz^2),\enskip \text{where $r=\sqrt{-1/\kappa}$.} \label{metric}\end{equation}
Clearly $H^3_\kappa$ convergences to the Euclidean $\oR^3$ on subsets on which $z$ remains bounded.
\par
A point $(x,y,z)\in H_\kappa$ corresponds to the curve in $\oP^1\times \oP^1-\ol\Delta$, given by the equation
\begin{equation} \frac{x+iy}{r} + \frac{x^2+y^2+(z+r)^2}{r^2}\zeta - \eta -\frac{x-iy}{r}\zeta\eta=0 \label{curve}\end{equation}
(see \cite{MS} or \cite{NR} for the case $\kappa=1$).  Multiplying \eqref{curve} by $r$ and regrouping yields
\begin{equation}
 (x+iy) +2z\zeta-(x-iy)\zeta\eta+r(\zeta-\eta)+\frac{x^2+y^2+z^2}{r}\zeta=0.\label{curve2}
\end{equation}
We view this curve as embedded in $\sX_{1/r}$. 
If we change coordinates back to $\zeta,w$, where $\eta=\zeta+tw=\zeta+w/r$, and let $r\to \infty$, we obtain the limit curve $ (x+iy) +2z\zeta-(x-iy)\zeta^2-w=0$, which corresponds to $(x,y,z)\in \oR^3$. 
\par
Thus, $\sigma$-invariant $(1,1)$-curves in $\sX_t$ correspond to points of  $H^3_{t^2}$, with the correspondence given by \eqref{curve} or \eqref{curve2}, with $r=t^{-1}$. Moreover:
\begin{proposition}
 Let $C_n\subset \sX_{t_n}$ be a sequence of  $\sigma$-invariant $(1,1)$-curves, given by the equation  $a^{00}_n+a^{10}_n\zeta-\frac{1}{t_n}
 \eta+a^{11}_n\zeta\eta=0$. Suppose that $t_n\to 0$, and that $a^{00}_n,a^{11}_n, a^{10}_n-1/t_n$ have finite limits $A^{00},A^{11},A^{10}$. Then the points $(x_n,y_n,z_n)\in H^3_{1/t_n^2}$ corresponding to $C_n$ converge to the point $(x,y,z)\in \oR^3$ corresponding to the limit curve $w=A^{00}+A^{10}\zeta+A^{11}\zeta^2$ in $\sX_0$.\hfill $\Box$
\end{proposition}

\section{Convergence of monopoles\label{conv1}}

We recall \cite{A, Hit1, MS} that Euclidean or hyperbolic monopoles correspond to algebraic curves in corresponding minitwistor spaces, i.e. in $T\oP^1$ or in $\oP^1\times \oP^1-\ol \Delta$. These curves satisfy certain conditions, the most important of which is the triviality of certain line bundle. In the case of hyperbolic monopoles of charge $k$ and mass $m$, this is the line bundle $\sO(2m+k,-2m-k)\simeq \sO(1,-1)^{2m+k}$ on $\oP^1\times \oP^1-\ol \Delta$ (the line bundle $\sO(1,-1)$ is topologically trivial on $\oP^1\times \oP^1-\ol \Delta$ and can be raised to an arbitrary complex power). In the case of Euclidean monopoles of charge $k$, this is the line bundle $L^{2}$, with transition function $\exp(-2w/\zeta)$ in the standard coordinates $\zeta$, $w\frac{\partial}{\partial\zeta}$.

We now combine these  line bundles to define certain line bundles on $\bar\sX$ (a partial compactification of $\sX$, described in Remark \ref{cpct}):
\begin{itemize}
 \item For every $a,b\in \oZ$, the line bundle $\sO(a,b)$ has transition function $\zeta^a\eta^b$ 
if $t\neq 0$ and $\zeta^{a+b}$ if $t=0$.
\item The line bundle $\sL$ has transition function  $\left(\frac{\eta}{\zeta}\right)^{-1/t}=\left(1+t\frac{w}{\zeta}\right)^{-\frac{1}{t}}$ if $t\neq 0$,  and $\exp(-w/\zeta)$ if $t=0$.
\end{itemize}
This last formula requires explanation, as $\left({\eta}/{\zeta}\right)^{-1/t}$ is not well defined unless $1/t$ is an integer. What the formula means is that $\sL|_{\sX_t}= \sO(1/t,-1/t)$.
\par
Thus, the line bundle $\sL^{2}(k,-k)\simeq \sL^{2}\otimes \sO(k,-k)$ restricted to $\sX_t$ is precisely the line bundle, the triviality of which is required for spectral curves of monopoles (this remains true even for massless monopoles, i.e. when $t=\infty$). 

\begin{remark} Once again, the line bundle $\sL$ is essentially given in \cite{A}, although only  its restrictions  $\sL|_{\sX_t}$ are discussed there.
\end{remark}

\medskip

We now wish to relate convergence of hyperbolic monopoles to a Euclidean one to convergence of corresponding spectral curves.

\subsection{Families of monopoles\label{fm}}
To capture the limiting behaviour of hyperbolic monopoles as their mass tends to infinity, we use a parameterised version of the Atiyah-Hitchin-Ward correspondence \cite{AW,Hit1}. We consider the following manifold:
\begin{equation}
 \sH=\{(x,y,z,t)\in \oR^4;\; zt>-1\},\label{sH}
\end{equation}
together with the projection $p:(x,y,z,t)\mapsto t$. We equip the vertical bundle $\Ker dp$ with the metric \eqref{metric}, where $r=1/t$. Thus, each fibre $\sH_t=p^{-1}(t)$ is has constant curvature $-t^2$. We shall consider vertical parts of other differential-geometric objects, e.g. vertical $p$-forms $\Omega^p_V\sH$ and vertical forms with values in a vector bundle $\Omega^p_V(E)=\Gamma\left(\Lambda^p_V\sH\otimes E\right)$.  
\par
We observe that, since each $\sH_t$, $t\neq 0$, is an upper half-space, we have chosen a preferred point on the ideal boundary of each hyperbolic space. Moreover, these points converge to a preferred point on the sphere at infinity of $\oR^3\simeq \sH_0$. It is given by the positive $z$-direction.
\par
We also observe that the ideal boundary of each $\sH_t$ varies smoothly with $t$, and so, we can add the sphere at infinity $S^2_\infty(t)$ for each $t$ to obtain a smooth manifold with boundary, which we denote by $\bar{\sH}$.
\par
We can now consider the Atiyah-Hitchin-Ward correspondence in families. On the one hand, a family of monopoles on hyperbolic spaces of curvature $-t^2$ (Euclidean for $t=0$), of some differentiability class $C^r$ in $t$ (we include here $C^\omega$, i.e. analytic in $t$), can be viewed as an object defined on a principal $SU(2)$-bundle $P$ over $\sH$: a pair $(\nabla,\Phi)$ consisting of a section $\Phi$ of $\ad P$ and  partial connection $\nabla$ (if $E$ is the associated rank $2$ complex vector bundle, then $\nabla:\Gamma(E)\to \Omega_V^1(E)=\Gamma\left(T_V^\ast \sH \otimes E\right)$), which define a monopole on each $\sH_t$. The $C^r$-dependence on $t$ translates into $\nabla$ and $\Phi$ being $C^r$ with respect to $t$ for $t\geq 0$.
\par
Moreover, since the Atiyah-Ward correspondence involves scattering along geodesics (\cite{Hit1,MS}), we must require that the boundary conditions on the sphere at infinity $S^2_\infty(t)$ also vary $C^r$ with $t$. In other words, $\Phi|_{S^2_\infty(t)}$ and the $U(1)$-connection $\nabla^0$ on $S^2_\infty(t)$ - the boundary value of a hyperbolic or Euclidean monopole (cf. \cite[\S5]{BA}) - must be $C^r$ in $t$.   

\subsection{Families of spectral curves\label{fam}}
  
The Atiyah-Hitchin-Ward correspondence associates to a monopole on $\sH_t$ (which is $ H^3_{t^2}$ if $t\neq 0$, and $\oR^3$ if $t=0$) a holomorphic vector bundle $\tilde E$ on $\sX_t$, trivial on curves of the form \eqref{curve2}, fitting into the following two exact sequences (over $\sX_t$):
$$ 0\to \sL(0,-k)\to \tilde E\to \sL^\ast(0,k)\to 0,$$
$$ 0\to \sL^\ast(-k,0)\to \tilde E\to \sL(k,0)\to 0.$$
The spectral curve $S$ of the monopole is defined a subscheme where the two extensions coincide. Clearly, the line bundle $\sL^2(k,-k)$ is trivial on $S$. Furthermore, the (real) section of $\sL^2(k,-k)|_S$ determines the above extensions, and consequently the monopole.
\par
Going through the proofs in \cite{Hit1} and \cite{MS} shows that these constructions depend smoothly on $t$. Thus a $C^r$-family of framed monopoles with $C^r$ boundary conditions, as described in the previous subsection, will produce a family of spectral curves $S_t\subset \sX_t$ together with section $s(t)$ of $\sL^2(k,-k)$ also $C^r$ in $t$. The inverse construction works as well. We can state the correspondence in a particularly elegant way, if we employ Nash's idea \cite{Nash}, of using the sections $s(t)$ to lift spectral curves to the total space of $\sL^2(k,-k)$.  To this end we shall denote by $\sZ_k$ the total space of $\sL^2(k,-k)$ over $\bar\sX$ without the zero section, and by $p:\sZ_k\to \oP^1$ the composition of the projection $\sZ_k\to \bar \sX$ and $\pi:\bar\sX\to \oP^1$. We denote by $\tau$ the real structure on $\sZ_k$ induced by $\sigma$ (which is defined on $\bar\sX$).
\par
We write $\sZ_k^+$ for $p^{-1}\{t\geq 0\}$.
\begin{proposition}
 The Atiyah-Hitchin-Ward correspondence induces a 1-1 correspondence between
\begin{itemize}
 \item[(i)] $C^r$ (in $t$) families of framed monopoles, as described in the previous subsection, and 
\item[(ii)] $\tau$-invariant surfaces $\sS$ in  $\sZ_k^+$, $C^{r}$ in the $t$-variable, and such that $\sS_t=\sS\cap p^{-1}(t)$ is a lift of a monopole spectral curve (i.e. it satisfies the additional conditions required of a curve in $\sX_t$ to be the spectral curve of a monopole).
\end{itemize}   
\end{proposition}

\section{Convergence of pluricomplex structures\label{cp}} 
\subsection{Pluricomplex structures\label{plu}} We briefly recall the definition and basic properties of a pluricomplex structure \cite{BS}.
\par
Let $V$ be a $2n$-dimensional real vector space. Associated to $V$ is its twistor space $\sJ(V)\simeq GL(2n,\oR)/GL(n,\cx)$, i.e. the space of complex structures on $V$. We denote by  $\sV^{1,0}$   the tautological $i$-eigenbundle over $\sJ(V)$. A {\em pluricomplex structure} is an embedding $K:\oP^1\to\sJ(V)$ such that both $K^\ast \sV^{1,0}$  and the quotient bundle $\cx^{2n}/K^\ast \sV^{1,0}$ split into line bundles of degree $-1$. This forces $n$ to be even.\par
To a pluricomplex structure one can associate a $\sigma$-sheaf $\sF$ on $\oP^1\times \oP^1-\ol\Delta$ 
 by considering a second pluricomplex structure $\wh K=-K\circ\sigma$ and the embedding 
$$ (K\circ \pi_1)^\ast  \sV^{1,0}\oplus (\wh K\circ \pi_2)^\ast  \sV^{1,0}\hookrightarrow \cx^{2n},$$
of sheaves on $\oP^1\times \oP^1$, where $\pi_1,\pi_2:\oP^1\times \oP^1\to \oP^1$ are the two projections. The sheaf $\sF$ is defined as the cokernel of this map and called the {\em characteristic sheaf} of the pluricomplex structure $\sF$. Its (scheme-theoretic) support $S$ is a compact curve in $\oP^1\times \oP^1$, called the {\em characteristic curve} of $K$. A pluricomplex structure is hypercomplex if and only if $S$ is the diagonal in $\oP^1\times \oP^1$ (i.e. $(S,\sO_S)\simeq (\Delta,\sO_\Delta)$). Since both $K^\ast \sV^{1,0}$ and  $\wh{K}^\ast \sV^{1,0}$ split into $\sO(-1)$'s, we can rewrite the exact sequence defining $\sF$ as
\begin{equation}0\to \sO(-1,0)^{\oplus n}\oplus \sO(0,-1)^{\oplus n}\stackrel{M(\zeta,\eta)}{\longrightarrow}\sO^{\oplus 2n}\to \sF\to 0. \label{F}\end{equation}
As shown in \cite{BS}, the map $M(\zeta,\eta)$ can be written (after fixing one complex structure) in the form
\begin{equation}
 M(\zeta,\eta)=\begin{pmatrix} X+\zeta Y & -1\\ \zeta & \eta\bar X-\bar Y\end{pmatrix},\label{M}
\end{equation}
and isomorphic pluricomplex structures correspond to orbits of such $M$ under the action of $GL(n,\cx)$ given by
$$ g. M=\begin{pmatrix} g & 0\\0 & \bar g\end{pmatrix} M \begin{pmatrix} {\bar g}^{-1} & 0\\0 & g^{-1}\end{pmatrix}.$$
The characteristic curve $S$ is now defined (as set) by the equation
\begin{equation}
S= \{(\zeta,\eta);\; \det M(\zeta,\eta)=0\}=\{(\zeta,\eta);\; \det\bigl((\eta\bar X-\bar Y)(X+\zeta Y)+\zeta I\bigr)=0\}.\label{S}
\end{equation}
In this setting, $K$ is hypercomplex if and only if $Y=0$ and $\bar X X=-1$.

\subsection{Limiting behaviour} Suppose that we have a sequence of pluricomplex structures on $V$ converging to a hypercomplex structure. The characteristic curves will converge to the diagonal $\Delta$ in $\oP^1\times \oP^1$, but there is additional data: if the tangent directions have a limit, then they will define a curve in $T\Delta\simeq T\oP^1$. Similarly, we expect that the characteristic sheaves $\sF$ will converge to a $1$-dimensional sheaf on $T\Delta$. In other words, we expect that a hypercomplex structure, which arises as a nontrivial limit of pluricomplex structures, carries more information than just the hypercomplex structure itself.
\par
Let us be more precise.
 Let $K_t$, $t\in I$, be a family of pluricomplex structures $K_t$, where $0\in I\subset \oR$  and $0$ is an accumulation point of $I$. We suppose  that $K_0$ is a hypercomplex structure and $K_t\to K_0$ as $t\to 0$. In addition, we suppose that $(K_t-K_0)/t$ has a limit $A:\oP^1\to \End V$ as $t\to 0$. Observe that $A(\zeta)\in T_{K_0(\zeta)}\sJ(V)$, i.e. $A(\zeta)$ anti-commutes with $K_0(\zeta)$ for every $\zeta$. 

Using the matrix description of the previous subsection, let us identify $K_t$ with pairs of matrices $(X_t,Y_t)$. Thus $Y_0=0$ and $\bar{X}_0X_0=-1$. We write $P=\lim_{t\to 0}\frac{X_t-X_0}{t}$, $Q=\lim_{t\to 0}\frac{Y_t}{t}$. Modulo $t^2$, equation \eqref{S} defining the characteristic curve can be rewritten as 
$$ \det\left( (\zeta-\eta) +t\eta \bar X_0(P+\zeta Q)+t(\eta\bar P-\bar Q)X_0\right)=0.$$
Replacing $\eta-\zeta=tw$, and dividing the above by $t$, we see that the characteristic curves $S_t$ converge to an $S_0\subset T\oP^1$ defined by
$$ S_0=\left\{(\zeta,w)\in T\oP^1;\, \det\bigl(wI + \bar Q X_0 -(\bar X_0 P+\bar P X_0)\zeta -\bar X_0Q\zeta^2\bigr)=0\right\}.$$
In other words, the characteristic curves define a ``surface'' in $\sX_I=\pi^{-1}(I)\subset \sX$, which is $C^1$ at $t=0$. To see the corresponding convergence of sheaves $\sF_t$, we observe that the characteristic sheaf of a pluricomplex structure $K$ is also the cokernel of the map
$$ (K\circ \pi_1)^\ast  \sV^{1,0}\rightarrow \cx^{2n}/(\wh K\circ \pi_2)^\ast,$$
defined as the composition of  the embedding into $\cx^{2n}$ and the natural projection. If we rewrite this map in terms of $X$ and $Y$, we  shall obtain:
\begin{equation}
 0\to \sO(-1,0)^{\oplus n}\stackrel{C(\zeta,\eta)}{\longrightarrow}\sO(0,1)^{\oplus n}\to \sF\to 0,\label{F2}
\end{equation}
where $C(\zeta,\eta)=(\eta\bar X-\bar Y)(X+\zeta Y)+\zeta I$. Thus, for $X_t,Y_t$ as above,
\begin{equation} \lim_{t\to 0}\frac{C_t(\zeta,\eta)}{t}=-wI - \bar Q X_0 +(\bar X_0 P+\bar P X_0)\zeta +\bar X_0Q\zeta^2,\label{C}\end{equation}
and the characteristic sheaves $\sF_t$ converge to a sheaf $\sF_0$ on $T\oP^1$ defined as cokernel of this map. 
\par
We wish to formalise the above statements using the manifold $\sX$ defined in section \ref{sX}. In order to avoid talking about ``sheaves which are $C^1$ in $t$-variable'', let us assume that $I$ is an open interval and that $K_t$ defines an analytic map $K:I\times \oP^1\to\sJ(V)$. We can view $K$ and other objects defined on $I$  as living in the CR-category. 
We can the rephrase the above as follows.
\begin{proposition} Let $I\subset \oR$ be an open interval containing $0$. There is 1-1 correspondence between
\begin{itemize}
 \item[(i)] CR-maps $K:I\times \oP^1\to\sJ(V)$ such that $K(0,\cdot)$ is hypercomplex and $K(t,\cdot)$ is pluricomplex for $t\neq 0$, and
\item[(ii)] coherent CR-$\sigma$-sheaves $\sF$ on $\sX_I=\pi^{-1}(I)\subset \sX$ such that, for every $t\in I$, $\sF|_{\sX_t}$ is $1$-dimensional and $$H^\ast(\sX_t,\sF(-2,0))=H^\ast(\sX_t,\sF(0,-2))=H^\ast(\sX_t,\sF(-1,-1))=0,$$
 where $\sF(a,b)=\sF\otimes \sO(a,b)$, and $\sO(a,b)$ has been defined at the beginning of \S\ref{conv1}. \hfill $\Box$
\end{itemize}\label{families}
\end{proposition}

\subsection{The geometry of pluricomplex limits\label{LHC}} As already mentioned above, a hypercomplex structure, which arises as a $C^1$-limit of pluricomplex structure carries more information than just a triple of anticommuting structures. The limit $A :\oP^1\to \End V$ of $(K_t-K_0)/t$ is a map from $\oP^1\to T\sJ(V)$, such that $A(\zeta)\in T_{K_0(\zeta)}\sJ(V)$, i.e. $A(\zeta)$ anti-commutes with $K(\zeta)$ for every $\zeta$. Let us therefore adopt the following definition:
\begin{definition} Let $V$ be a vector space. An l-hypercomplex structure is a holomorphic map $A:\oP^1\to T\sJ(V)$ such that $\pi\circ A$ is a hypercomplex structure (here $\pi:T\sJ(V)\to\sJ(V)$ is the base map).\label{l}
\end{definition}
In other words, an l-hypercomplex structure (``l'' for limit) is a lift of a hypercomplex structure to $T\sJ(V)$. Clearly, any hypercomplex structure defines trivially an l-hypercomplex structure (as a zero section), but the previous section shows that a limit of pluricomplex structures will usually be non-trivial.
\par
Let us consider such an $l$-hypercomplex structure $A$, and denote by $K$ the hypercomplex structure $\pi\circ A$. We extend linearly $A(\zeta)$ to $V^\cx$ and consider the map $\tilde A(\zeta):V^{1,0}_\zeta\to V^\cx/V^{1,0}_\zeta$, obtained as restriction of $A(\zeta)$ composed with the natural projection ($V^{1,0}_\zeta$ consists of vectors of type $(1,0)$ for $K(\zeta)$). Thus, we obtain a bundle map $\tilde A:\sV^{1,0}\to V^\cx/\sV^{1,0}$ over $K(\oP^1)$.
\par
In addition, if $L_\zeta=\langle K(\zeta^\prime)\in \End V;\, K(\zeta^\prime)\circ K(\zeta)=-K(\zeta)\circ K(\zeta^\prime)\rangle$, then $L$ is a holomorphic line bundle over $K(\oP^1)$, isomorphic to $T\oP^1$. We pull-back $\sV^{1,0}$ and $V^\cx/\sV^{1,0}$ to the total space of $L$. Since any $K(\zeta^\prime)$, which anticommutes with $K(\zeta)$, maps $V^{1,0}_\zeta$ to $V^{0,1}_\zeta$, we obtain a tautological map $\sI:\sV^{1,0}\to V^\cx/\sV^{1,0}$ on $\rm{Tot}(L)\simeq T\oP^1$. The map $\tilde A-\sI$ on $T\oP^1$ is an isomorphism for a generic fibre, and we obtain an exact sequence of sheaves:
\begin{equation}
 0\to V^{1,0}\stackrel{\tilde A-\sI}{\longrightarrow} V^\cx/V^{1,0}\to \sF\to 0.\label{lF} 
\end{equation}
If $K$ is written as in \eqref{M}, using matrices $X,Y$, so that $A(\zeta)$ is identified with $P+\zeta Q$ for some $n\times n$ complex matrices $P,Q$, then $\tilde A$ has the form as in \eqref{C}:
$$ \tilde{A}= -\bar Q X +(\bar X P+\bar P X)\zeta +\bar XQ\zeta^2,
$$
and $\sI$ is simply $w\rm{I}$, where $w$ is the fibre coordinate of $T\oP^1$. Thus \eqref{lF} can be rewritten as an exact sequence of sheaves on $T\oP^1$:
\begin{equation}
 0\to \sO(-1)^{\oplus n}\stackrel{\tilde A(\zeta)-w\rm{I}}{\longrightarrow}\sO(1)^{\oplus n}\to \sF\to 0.\label{lF2} 
\end{equation}
Observe that $\tilde A$ is compatible with the natural real structure, i.e. $\tilde A(-1/\bar\zeta)=-\ol{\tilde A(\zeta)}/\bar\zeta^2$. 
We can show, similarly to \cite{BS}:
\begin{proposition}
 There is a 1-1 correspondence between l-hypercomplex structures on $\cx^{2n}$ and coherent $1$-dimensional $\sigma$-sheaves $\sF$ on $T\oP^1$ such
 that $h^0(\sF)=2n$ and $H^\ast(\sF(-2))=0$.\hfill $\Box$
\end{proposition}
As for pluricomplex structures, we call $\sF$ the {\em characteristic sheaf} of an l-hypercomplex structure $A$ and its support $S$ the {\em characteristic curve} of $A $. The degree of $A$ is the degree of $S$.
\par
L-hypercomplex structures of degree $1$ are simply hypercomplex structures; for higher degrees, however, they carry additional information.
\par
We can recover $A(\zeta)$ from $\sF$ as follows. The vector space $V$ is identified with $\sigma$-invariant sections of $\sF$ and the complex structure $J_{\zeta_0}$ corresponding to $\zeta_0\in \oP^1$ is obtained via the isomorphism (evaluation) $H^0(\sF)^\sigma\simeq H^0(D_{\zeta_0},\sF)$, where $D_{\zeta_0}$ is the divisor cut out  by $\zeta-\zeta_0$. Assume, for simplicity, that $D_{\zeta_0}\cap S$ consists of distinct points $u_1,\dots,u_k$ in $T_{\zeta_0}\oP^1$. According to the above discussion, we can identify $T_{\zeta_0}\oP^1$ with the complex line generated by complex structures $J_\zeta$ anticommuting with $J_{\zeta_0}$. If $J^\prime$ is one of those, then we can identify the points $u_1,\dots,u_k$ with $w_1J^\prime,\dots,w_kJ^\prime$, for some complex numbers $w_1,\dots,w_k$. We have:
$$ V\simeq  H^0(D_{\zeta_0},\sF)\simeq \bigoplus_{i=1}^k  H^0(\{u_i\},\sF).$$
Let $\pi_i$ denote the projection onto the $i$-th direct summand $ H^0(\{u_i\},\sF)$. Observe that each of these summands is $J_{\zeta_0}$-invariant, and, hence, can be thought of as a complex subspace (for $J_{\zeta_0}$). Then, for $v\in V$:
\begin{equation} 
 A(\zeta_0)(v)=\sum_{i=1}^k w_i\pi_i(J^\prime v).\label{AAA}
\end{equation}

\subsection{Strongly integrable l-hypercomplex structures} We recall, from \cite{BS}, that a pluricomplex manifold $M$ is called {\em strongly integrable} if $M$ can be recovered as the parameter space of curves in a complex manifold $Z$ fibreing over $\oP^1\times \oP^1-\ol\Delta$, in such a way that the curve $C_m$ corresponding to $m\in M$ is a lift of the characteristic curve $S_m\subset \oP^1\times \oP^1-\ol\Delta$ of the pluricomplex structure on $T_mM$, and  the normal bundle of $C_m$ in $Z$ is isomorphic to the characteristic sheaf $\sF_m$.
\par
Similarly to the construction of $Z$ in \cite{BS}, consider, for an l-hypercomplex manifold $M$, the following fibration over $M$:
\begin{equation} Y=\{(m,x)\in M\times T\oP^1;\enskip x\in S_m\}\stackrel{\nu}\longrightarrow M ,\label{Y}\end{equation}
 and the corresponding map $p:Y\to T\oP^1$. We shall assume that $Y$ is smooth (it is certainly smooth whenever  $S_m$ is smooth) and that $p$ is a submersion onto its image. For every $\zeta\in \oP^1$, we consider the subbundle $T_\zeta^{1,0}M$ of vectors of type $(1,0)$ for the complex structure $K(\zeta)$, and the bundle map $\tilde A(\zeta): T_\zeta^{1,0}M\to T^\cx M/T_\zeta^{1,0}M$.
\begin{definition} An l-hypercomplex structure $A$ is said to be {\em strongly integrable} if 
$T^{1,0}_\zeta+\im(w-\tilde A(\zeta))$ is a subbundle of  $T^\cx p^{-1}(\zeta,w)$ for all $(\zeta,w)\in \im p;$ and the resulting distribution $\sD\subset T^\cx Y$ is involutive (i.e. $[\sD,\sD]\subset \sD$). A manifold equipped with a strongly integrable l-hypercomplex structure will be called {\em strongly l-hypercomplex}.\end{definition}
\begin{definition} The twistor space $Z$ of a strongly integrable l-hypercomplex structure is the space of leaves $Y/\sD^\oR$, where  $\sD^\oR=\sD\cap \overline{\sD}$.
\end{definition}
$Z$ is a holomorphic foliation, with a holomorphic map $\rho:Z\to TP^1$, induced by $p$. Moreover, $Z$ has an antiholomorphic involution $\tau$, which covers the involution $\sigma$, $\sigma(\zeta,w)=(-1/\bar\zeta, -\bar w/\zeta^2)$. Each characteristic curve $S_m$, $m\in M$, descends to $Z$ and its normal bundle in $Z$ is isomorphic to $\sF_m$.
\par
As for pluricomplex structures \cite[Theorem 6.6]{BS}, we have:
\begin{theorem} Let $Z$ be a complex manifold equipped with  a holomorphic map $\rho:Z\to T \oP^1$ and an antiholomorphic involution $\tau:Z\to Z$ such that $\rho\circ\tau=\sigma\circ\rho$. Let $\sC$ be a maximal connected subset of the Douady space $H_Z$, consisting of complex subspaces  $S$ satisfying the following conditions:
\begin{itemize} 
\item[(i)] 
$\rho_{|S}$ is a biholomorphism onto its image $\rho(S)$, which is a compact divisor of $\sO(2k)$ for some $k\in \oN$;
\item [(ii)] the sheaf $\sF=(\rho_{|S}^{-1})^\ast \sN_{S/Z}$ satisfies
$H^\ast\bigl(\rho(S), \sF(-2)\bigr)=0.$
\end{itemize}
 Then $\sC$ is a smooth manifold and its $\tau$-invariant part $M=\sC^\tau$   is equipped with a natural strongly integrable l-hypercomplex structure. Moreover the twistor space of $M$ is the corresponding open subset of $Z$.\hfill $\Box$\label{twistor}
\end{theorem}

\subsection{Strongly integrable deformations}

We can consider strongly integrable pluricomplex deformations of strongly integrable l-hypercomplex manifolds, similar to the linear deformation described in Proposition \ref{families}. In order to state the precise result, we need a definition.
\begin{definition} A {\em strongly integrable pluricomplex deformation} is an analytic manifold $\sM$, together with the following data:
\begin{itemize}
 \item An analytic submersion $p:\sM\to I$, where $I$ is an interval containing $0$;
\item A CR-map $K:\oP^1\to \sJ(\Ker dp)$, such that $K$ induces a strongly integrable pluricomplex structure on $\sM_t=p^{-1}(t)$ for $t\neq 0$, and a hypercomplex structure on $\sM_0=p^{-1}(0)$. Moreover, $\frac{dK}{dt}|_{t=0}$ is a  strongly integrable l-hypercomplex structure on $\sM_0$.
\end{itemize}
\end{definition}
The following result follows from the definition and basic deformation theory.
\begin{theorem} There is a natural 1-1 correspondence between strongly integrable pluricomplex deformations $p:\sM\to I$, {\em and} CR-submersions $\pi:\sZ\to I$, equipped with a CR-map $\rho:\sZ\to \sX_I=\pi^{-1}(I)\subset \sX$, and an antiholomorphic involution $\tau$ such that $\rho\circ\tau=\sigma\circ\rho$. \hfill $\Box$\label{sZ}
\end{theorem}

\begin{remark}
 We can always extend everything into the complex domain, and so we can consider instead complex manifolds $\sZ$ with a holomorphic submersion $\sZ\to D$, where $D$ is a neighborhood of $I$ in $\cx$. 
\end{remark}

\subsection{Geometry of monopoles revisited}

Moduli spaces of monopoles fit very well into the above theory. The moduli space $M_k$ of framed Euclidean $SU(2)$-monopoles of charge $k$ has a strongly integrable l-hypercomplex structure. Its twistor space $Z_k$ is the total space of the line bundle $L^2$ over $T\oP^1$ without the zero section. The moduli space $M_k$ is obtained as in Theorem \ref{twistor}, where the complex subspaces $S$ are lifts of spectral curves of monopoles to $Z_k$. Results of Nash \cite{Nash} show that the conditions of this Theorem (in particular, the condition (ii) on the normal bundle) are satisfied.
\par
The resulting  l-hypercomplex structure on $M_k$ is essentially described in \eqref{AAA}. Fix $\zeta_0$ and identify the moduli space with based rational maps $p(z)/q(z)$ of degree $k$. On the subset where the poles $\eta_1,\dots,\eta_k$ are distinct, we have local coordinates $\eta_i,p(\eta_i)$, $i=1,\dots,k$. The projection $\pi_i$ used in \eqref{AAA} is simply the map $p(z)/q(z)\mapsto (\eta_i,p(\eta_i))$.
\par
Furthermore, moduli spaces of hyperbolic monopoles with varying mass form a strongly integrable deformation of $M_k$, as defined in the previous subsection. On the one hand, we can, as in \S\ref{fm}, define a smooth manifold $\sM_k$ fibreing over $\oR$, such that the fibre over $t$ is the moduli space of framed monopoles on the hyperbolic $3$-space with curvature $-1/t^2$. On the other hand, we have the complex  manifold
$\sZ_k$, defined (\S\ref{fam}) as the total space of $\sL(k,-k)$ over $\bar\sX$ without the zero section, equipped with a fibration $p:\sZ_k\to \oP^1$.
We would like to say that $\sM_k$ and $\sZ_k$ are an example of $\sM$ and $\sZ$ in Theorem \ref{sZ}, but at present we can only say this about their open subsets. The problem is that we do not know whether the pluricomplex structure on moduli spaces of hyperbolic monopoles is everywhere nonsingular. More precisely, we do not know whether the normal bundle $N$ in $p^{-1}(t)$ of a lift $S$ of a spectral curve of a hyperbolic monopole always satisfies $H^\ast(S,N(-2,0))=0$ (cf. \cite{BS} for partial results and discussion). Still,  since  the normal bundle of a lift of the spectral curve of a Euclidean monopole to $p^{-1}(0)$ satisfies $H^\ast(S,N(-2))=0$ \cite{Nash}, it follows that there is an open neighbourhood of $p^{-1}(0)$ in $\sZ_k$ where lifts of hyperbolic spectral curves satisfy $H^\ast(S,N(-2,0))=0$, i.e. the geometry is pluricomplex. Thus, there exist open neighbourhoods $\sM_k^\ast,\sZ_k^\ast$ of the fibres over $t=0$ in $\sM_k$ and in $\sZ_k$ (i.e. of $M_k$ and $Z_k$) which correspond to $\sM$ and $\sZ$ in Theorem \ref{sZ}.

\begin{remark}
 We finish with few words about the limit $t\to \infty$, i.e. when the mass of a monopole tends to zero. It has been shown in \cite{BS} that the moduli space of massless monopoles of charge $k$ is not pluricomplex. What happens is that the matrices $X$ and $Y$ defined in \eqref{M} acquire singularities as $t\to \infty$. The limiting geometry is weaker than pluricomplex: in the notation of \S\ref{plu}, it is given by an embedding $K:\oP^1\to \sJ(V)$ such that $K^\ast\sV^{1,0}$ still splits into line bundles of degree $-1$, but the quotient bundle $\cx^{2n}/K^\ast \sV^{1,0}$ does not.
\end{remark}

\end{document}